\documentclass[11pt,a4paper]{amsart}
\usepackage{amssymb,amsmath}
\usepackage{cancel}
\usepackage{palatino}
\usepackage{longtable}
\usepackage[all]{xy}
\usepackage{tikz-cd}
\usepackage{color}
\numberwithin{equation}{section}
\newtheorem{theorem}{Theorem}[section]     
\newtheorem{definition}[theorem]{Definition}
\newtheorem{defi}[theorem]{Definition}

\newtheorem{lemma}[theorem]{Lemma}

\newtheorem{remark}[theorem]{Remark}

\def\dna{d_{\nabla}}
\def\dlna{d_{L\nabla}}
\def\na{\nabla}

\def\O{\mathcal{O}_M}
\def\T{\mathcal{T}_M}
\def\OO{\Omega_M}
\def\d{\partial}

\def\dna{d_{\nabla}}
\def\dlna{d_{L\nabla}}
\def\Ll{L_{\lambda}}
\def\na{\nabla}
\def\nas{\nabla^*}
\def\proof{\noindent\hspace{2em}{\itshape Proof: }}
\def\QEDclosed{\mbox{\rule[0pt]{1.3ex}{1.3ex}}} 

\def\QED{\QEDclosed} 
\def\endproof{\hspace*{\fill}~\QED\par\endtrivlist\unskip}
\newcommand{\eqa}{\begin{eqnarray}}
\newcommand{\eeqa}{\end{eqnarray}}
\newcommand{\beq}{\begin{equation}}
\newcommand{\eeq}{\end{equation}}

\numberwithin{equation}{section}

\begin{document}
\title{Bi-flat F-structures as differential bicomplexes and Gauss-Manin connections}

\author{Alessandro Arsie}
\address{A.~Arsie:\newline Department of Mathematics and Statistics, The University of Toledo,\newline 2801W. Bancroft St., 43606 Toledo, OH, USA}
\email{alessandro.arsie@utoledo.edu}

\author{Paolo Lorenzoni}
\address{P.~Lorenzoni:\newline Dipartimento di Matematica e Applicazioni, Universit\`a degli Studi di Milano-Bicocca, 
Via Roberto Cozzi 53, I-20125 Milano, Italy and INFN sezione di Milano-Bicocca}
\email{paolo.lorenzoni@unimib.it}

\date{}

\maketitle

\begin{abstract}
We show that a bi-flat  F-structure $(\nabla,\circ,e,\nabla^*,*,E)$ on a manifold $M$ defines a differential bicomplex  $(d_{\nabla},d_{E\circ\nabla^*})$  on forms  with value  on  the tangent sheaf of the manifold. Moreover,  the sequence of vector fields defined recursively by $d_{\nabla}X_{(\alpha+1)}=d_{L\nabla^*}X_{(\alpha)}$ coincide with the coefficients of the formal expansion of the flat local sections of a family of flat  connections  $\nabla^{GM}$ associated with the bi-flat structure.  In the case of Dubrovin-Frobenius manifold the connection $\nabla^{GM}$  (for  suitable  choice of  an auxiliary  parameter) can be identified with  the Levi-Civita connection of  the flat pencil  of  metrics defined by the invariant metric and the intesection form. 
\end{abstract}
\tableofcontents

\section{Introduction}
Bi-flat F-manifolds  $(\nabla,\circ,e,\nabla^*,*,E)$ have been introduced in \cite{ALjgp}  as a generalization of the notion of Dubrovin-Frobenius manifolds, initially motivated by the possibility to include in this set-up also non-Hamiltonian integrable PDEs. This has also led to the realization of such a structure on the orbit spaces of some well-generated complex reflection groups (\cite{KMS15,ALcomplex,KMS}). It  is worth mentioning also the relations with Painlevé trascendents and their generalizations \cite{ALjgp,Limrn,ALmulti,KMS15,KM}, the special class of bi-flat structures related to Lauricella functions \cite{ALmulti,LP} and the related quasi-homogeneous formal solutions of the oriented associativity equations \cite{C}. Recently, in \cite{ABLR1} semisimple bi-flat F-manifolds appeared as suitable genus 0 counterpart to Dubrovin-Frobenius manifolds in homogeneous F-Cohomological Field Theories, and through the theory developed there it was possible to prove the existence of dispersive integrable hierarchies, not necessarily Hamiltonian associated to a semisimple bi-flat F-manifold, see \cite{ABLR2}. 

 Dubrovin-Frobenius manifolds are bi-flat F-manifolds equipped with an invariant flat metric satisfying suitable compatibility conditions. They have been introduced by B. Dubrovin in the 90's to provide a differential geometric set-up for topological quantum field theories and quantum cohomology. Nowadays, they appear as important tools in many areas, including  singularity theory, Gromov-Witten invariants, integrable systems of Hamiltonian PDEs, mirror symmetry, orbit spaces of Coxeter groups, just to name a few. 

A key ingredient in the theory of Dubrovin-Frobenius manifolds is provided by a flat pencil of metrics $g_{\lambda}=g-\lambda\eta$ and the associated (extended) Gauss-Manin system for the flat coordinates of the pencil (see for instance \cite{du93} and \cite{du97}). This is essential not only for the Dubrovin-Frobenius manifold under investigation {\em per se}, but also because these flat coordinates are related to the Hamiltonian functionals of the associated bi-Hamiltonian hierarchy. These functionals are obtained by means of bi-Hamiltonian   recursion \cite{Magri}:
\begin{equation}\label{LMrr}
d_{P_\eta} H_{\alpha+1}=d_{P_g} H_{\alpha}
\end{equation}
where $P_\eta$ and $P_g$ are the differential geometric local Hamiltonian operators
  associated with the metrics $\eta$ and $g$, whereas $d_{P_g}$ and $d_{P_\eta}$ are the corresponding differentials in Poisson cohomology. To be more precise, the bi-Hamiltonian structure defined by  these operators can be interpreted as a bi-complex $(d_{P_\eta},d_{P_g})$ on the space of local multivectors on the loop space of the Dubrovin-Frobenius manifold. Indeed, the fact that $P_\eta$ and $P_g$ define a bi-Hamiltonian structure can be rephrased as
\[ d_{P_\eta}^2=0,\qquad d_{P_g}^2=0,\qquad d_{P_\eta}d_{P_g}+d_{P_g}d_{P_\eta}=0.\]
Furthermore, the flat pencil of metrics $g_{\lambda}$  and the associated Gauss-Manin system for its flat coordinates play an important role in the construction of Dubrovin's superpotential \cite{du93} and of the associated spectral curve in the theory of topological  recursion  \cite{DNOPS}.   
\newline
\newline 
While many aspects of Dubrovin-Frobenius manifolds have been generalized to the case of bi-flat F-manifolds, especially in the semisimple case, a glaring omission so far has been the construction of the analog of  Gauss-Manin system for a bi-flat F-manifold. In principle, this seems a daunting task since there are no significant flat metrics on a general bi-flat F-manifold. However, the driving idea motivating our  study is  the observation that the Lenard-Magri recursion relations \eqref{LMrr} for vector fields on the loop  space of the manifold, in the case of hamiltonian densities depending only on the coordinates  on  the target manifold (and not on their derivatives with respect to the local coordinate on the loop) reduce to Lenard-Magri recursion relations
 associated to a differential-bicomplex on the space of vector valued  differential forms. Furthermore, the definition of this differential bicomplex involves only the Levi-Civita connections of the metric $\eta$ and $g$ but not the metrics themselves. This crucial fact allows us to define this bicomplex in the broader context of bi-flat F-manifolds, where the the flat metrics $\eta$ and $g$ are replaced by two ``compatible'' flat connections  $\nabla$ and $\nabla^*$. We show that:
\begin{itemize}
\item  exactly as compatibility of the Hamiltonian  operators is related to the anticommutativity of the associated differentials, the compatibility of  the connections 
 $\nabla$ and $\nabla^*$ appearing in the definition of bi-flat F-manifolds is related 
   to the anticommutativity of the differentials which define the associated bicomplex 
\[(d_{\nabla}, d_{E\circ \nabla^*}, \T\otimes_{\O}\OO^{\bullet})\] on the locally free $\O$-module of forms with values in $\T$ (here $E\circ$  is the operator of multiplication by the Euler vector field). This differential bicomplex can be thought as a non-trivial lifting of the Fr{\"o}licher-Nijenuis bicomplex $(d, d_{E\circ}, \OO^{\bullet})$ to $\T\otimes_{\O}\OO^{\bullet}$;
\item exactly as the axioms of Dubrovin-Frobenius manifolds ensure that the pencil
 of metrics $g_{\lambda}$ is flat and thus the associated Gauss-Manin system is well defined, the axioms of  bi-flat F manifolds ensure that the one-parameter family  of  connections $\nabla^{GM}$  defined as 
 \begin{equation*}\label{GMnoindexIntro}
\nabla^{GM}_XY=\nabla^*_XY+\lambda(\nabla^*_{(E\circ-\lambda e\circ)^{-1}X}Y-\nabla_{(E\circ-\lambda e\circ)^{-1}X}Y).
\end{equation*}
is flat and torsionless for any $\lambda$.
\end{itemize}
It turns out that two structures introduced above are closely related to each  other:  the local sections of $\T$ that are covariantly constant with respect to $\nabla^{GM}$ can be expanded in $\lambda^{-1}$ and the coefficients $X_{\alpha}$ of this expansion satisfy the recurrence relations 
\begin{equation*}\label{LMrr2}
d_{\nabla}X_{\alpha+1}=d_{E\circ\nabla^*}X_{\alpha}.
\end{equation*}
\newline
\newline
 The organization of the paper follows the vertices and the arrows of  this triangular diagram:
\newline
\newline 
$\xymatrix
   { & (M, \circ, \nabla, e, \ast, \nabla^\ast, E) \ar@{<->}[dl] \ar@{<->}[dr] & \\
    (d_{\nabla}, d_{E\circ\nabla^*}, \T\otimes_{\O}\OO^{\bullet}) \ar@{<->}[rr] & & \nabla_{GM}
   }
$
\newline
\newline
In Section 2, we provide some background material on flat an bi-flat F-manifolds (upper vertex). In Section 3, we introduce the a  family of Gauss-Manin connections (right lower vertex), we show how its flatness and torsion freeness are related the axioms of bi-flat F-manifolds (arrows in the right side of the triangle) and we briefly show how to recover the standard construction in the case of Dubrovin-Frobenius manifolds. In Section 4, we introduce the aforementioned differential bicomplex (left lower vertex) and we show how the properties of its differentials are related the axioms of bi-flat F-manifolds (arrows in the left side of the triangle). Moreover, we construct the associated Lenard-Magri chains and we compare them with the flat sections of the Gauss-Manin connections.   
\newline 
\newline
We work indifferently in the category of $C^{\infty}$ manifolds, or $C^{\omega}$ real analytic manifolds, or complex (holomorphic) manifolds. In the latter case, whenever we talk about a metric we mean a complex-valued non-degenerate bilinear form, while in the real case this is just a pseudo-Riemannian metric. We denote with $\O$ the structure sheaf of $M$, with $\T$ the tangent sheaf of $M$, with $\T^*$ the cotangent sheaf and with $\OO^{\bullet}$ the graded sheaf of $\O$-modules given by forms (smooth, real analytic, or holomorphic) on $M$. Of course, the grading of $\OO^{\bullet}$ is given by  $\OO^{\bullet}=\oplus_{k=1}^n \OO^k$, where $\OO^k$ is the $\O$-module of $k$-forms and $n=\text{dim}(M)$. 
We denote with $\circ$ the commutative associative product on $\T$ defined in \ref{defFmani} and with $\cdot$ the composition of morphisms. All the sheaves involved are locally free finite rank sheaves of $\O$-module, so it is not restrictive to think of them as vector bundles if one desires so. 

\section{Flat and bi-flat F manifolds}
F-manifolds have been introduced by Hertling and Manin in \cite{HM}.
\begin{definition}\label{defFmani}
An \emph{F-manifold} is a manifold $M$ equipped with
\begin{itemize}
\item a morphism of sheaves $\circ:\T\times\T\rightarrow \T$  which is $\O$-bilinear, and it gives rise to a commutative associative product on $\T$. Moreover, it satisfies the following identity:
\begin{align}
&[X\circ Y,W\circ Z]-[X\circ Y, Z]\circ W-[X\circ Y, W]\circ Z\label{HMeq1free}\\
&-X\circ [Y, Z \circ W]+X\circ [Y, Z]\circ W +X\circ [Y, W]\circ Z\notag\\
&-Y\circ [X,Z\circ W]+Y\circ [X,Z]\circ W+Y\circ [X, W]\circ Z=0,\notag
\end{align}
for all local sections  $X,Y,W, Z$ of $\T$, where $[X,Y]$ is the Lie bracket.
\newline
\item A distinguished global section $e$ of $\T$ such that 
\[e\circ X=X\] 
for all local sections $X$ of $\T$. This acts like the unit of $\circ$.
\end{itemize}
\end{definition}
In any coordinate system, the morphism $\circ$ can be expressed using its structure functions $c^k_{ij}$ via the formula $\partial_i\circ\partial_j=c^k_{ij}\partial_k$, where $\{\partial_1, \dots, \partial_{n}\}$, with $n=\text{dim}(M)$ is a local basis of $\T$ associated to the chosen coordinate system. Commutativity of $\circ$ translates into the equation 
\begin{equation}\label{comm1.eq}
X\circ Y=Y\circ X \text{ for all local sections $X,Y\in \T$ or } \;\; c^k_{ij}=c^k_{ji},
\end{equation}
while the associativity of $\circ$ is equivalent to 
\begin{equation}\label{asso1.eq}
(X\circ Y)\circ Z=X\circ (Y \circ Z) 
\end{equation}
for all local sections $X,Y\in \T$ or 
\begin{equation}\label{asso1.eqBis}
c^i_{sj}c^s_{km}=c^i_{sk}c^s_{jm}.
\end{equation}
Of course, since $X\circ \in \T\otimes_{\O}\OO^1$, we can also view the morphism $\circ$ as $\circ: \T\rightarrow \T\otimes_{\O}\T^*$.

Usually F-manifolds are equipped with additional structures. For instance, in \cite{DS} motivated by Dubrovin's duality the authors introduced the notion of  eventual identity.
\begin{definition}\label{defFwithEvId}
An \emph{F-manifold with an eventual identity} is an  F-manifold $M$ equipped with a global section $E$ of $\T$ satisfying
\beq
\mathcal{L}_E(\circ)(X,Y)=[e,E]\circ X\circ Y,
\eeq
for all local section $X,Y$ of $\T$.
\end{definition}
The presence of an eventual identity $E$ allows one  to introduce a new morphism of sheaves $*:\T|_{U}\times \T|_{U}\rightarrow \T|_{U}$: \[X*Y:=(E\circ)^{-1}X\circ Y,\] 
where $X$ and $Y$ are arbitrary local sections of $\T|_{U}$. Indeed, in general $*$ is not defined everywhere on $M$ but on the open set $U$ where the endomorphism $E\circ:\T \rightarrow \T$ is invertible. 

In \cite{ALimrn} it has been proved that if $E$ is an eventual identity, then the Nijenhuis torsion of the sheaf morphism $L=E\,\circ: \T\rightarrow \T$, defined as
$$N_{L}(X,Y)=[E\circ X,E\circ Y]+E\circ E\circ [X,Y]-E\circ [X,E\circ Y]-E\circ [E\circ X,Y],$$
vanishes.

A special case of eventual identities is the case of Euler vector fields.  
\begin{definition}\label{defFwithE}
An \emph{F-manifold with Euler vector field} is an F-manifold $M$ equipped with a global section of $\T$ satisfying
$$[e,E]=e,\qquad  \mathcal{L}_E \circ=\circ.$$
\end{definition}

Other additional structures have been introduced motivated by the theory of integrable systems of hydrodynamic type \cite{LPR}. There the additional datum is a connection $\nabla$ satisfying suitable compatibility conditions coming from integrability
  of the associated integrable hierarchies. When $\nabla$ is flat this notion coincides with the notion of
 flat F-manifolds introduced by Manin in \cite{manin}.
\begin{defi}
A flat F-manifold is an  $F$-manifold $(M, \circ, e)$ equipped with a connection
 $\nabla:\T\rightarrow \T\otimes_{\O}\OO^1 $ satisfying the following axioms:
\begin{enumerate}
\item The one parameter family of connections
$$\nabla-\lambda\circ$$
is flat and torsionless for any $\lambda$. 
\newline
\item
The vector field $e$ is covariantly constant: $\nabla e=0$.
\end{enumerate}
\end{defi}
Bi-flat F-manifolds are manifolds equipped with two ``compatible'' flat structures or flat F-manifolds equipped with a \emph{linear} Euler vector field that is
\begin{equation}\label{eulerlinear.eq}
\nabla^2_{X,Y} E=0, 
\end{equation}
for any  local  sections $X,Y$ of $\T$  (this notion, in turn, is equivalent to the notion of Saito structure without metric \cite{Sa}). For the equivalence of these definitions we refer to  Theorem 4.4 in \cite{ALcomplex} for the semisimple  case and Lemmas 4.2 and 4.3 in \cite{KMS} for the general case. The notion of bi-flat F-manifolds has been introduced
  in  \cite{ALjgp} motivated by  the study of twisted Lenard-Magri chains started in \cite{ALimrn}. 

\begin{defi}\label{biflatdefi}
A \emph{bi-flat}  F-manifold is a manifold $M$ equipped with a pair
 of connections $\nabla$ and $\nabla^{*}$ on $\T$, a pair of sheaves morphisms $\circ$ and $*$ from $\T\times \T\rightarrow \T$ inducing commutative and associative products on $\T$  and  a pair of global sections of $\T$, $e$ and  $E$, s.t.:
\begin{itemize}
\item $(\nabla,\circ,e)$ defines a flat F-structure on $M$.
\newline
\item $(\nabla^{*},*,E)$ defines a flat F-structure on (an open dense subset $U$ of) $M$. 
\newline
\item The two structures are related by the following conditions:
\newline
\begin{itemize}
\item $E$ is an Euler vector field for the first structure and at a generic point the operator $E\circ$ is assumed to be invertible. 
\newline
\item $*$ is the product defined by $E$, i.e. $X*Y=(E\circ)^{-1} X\circ Y.$
\newline
\item $(d_{\nabla}-d_{\nabla^{*}})(X\,\circ)=0,$ where the morphism $\dna:\T\otimes_{\O}\OO^k\rightarrow  T\otimes_{\O}\OO^{k+1}$ is usually called exterior covariant derivative of vector-valued differential forms and is defined via:
$$(\dna \omega)(X_0, \dots, X_k)=\sum_{i=0}^k (-1)^i \na_{X_i}(\omega(X_0, \dots, \hat{X}_i, \dots, X_k))+$$
$$+\sum_{0\leq i<j\leq k}(-1)^{i+j}\omega([X_i, X_j], X_0, \dots, \hat{X}_i, \dots, \hat{X}_j, \dots X_k).$$
\end{itemize}
\end{itemize}
\end{defi}
In the previous definition, $\nabla$ and $\circ$ are usually called the natural connection and the natural product respectively, while $\nabla^*$ and $*$ are usually named the dual connection and the dual product respectively. 

Let us remark that not all the axioms are independent. For instance the compatibility between the dual connection and the dual product follows from the other axioms \cite{ALmulti}. Furthermore, the dual connection is defined only at the points where the operator $E\circ$ is invertible. At these points the condition 
\beq\label{ahe1}
(d_{\nabla}-d_{\nabla^{*}})(X\,\circ)=0\qquad \text{for all local sections } X\in \T|_{U}
\eeq
 is equivalent to the condition
\beq\label{ahe2}
(d_{\nabla}-d_{\nabla^{*}})(X\,*)=0\qquad\  \text{for all local sections } X\in \T|_{U}.
\eeq
In local coordinates $\circ$ is represented by the tensor $c$ such that $\partial_j\circ\partial_k=c^i_{jk}\partial_i$ and analogously $*$ is represented by the tensor $c^{*}$ such that $\partial_j*\partial_k=c^{*i}_{jk}\partial_i$. Likewise we call $a^i_{jk}$ the Christoffel symbols of $\nabla$ and $b^i_{jk}$ the Christoffel symbols of $\nabla^*$. Then the two relations \eqref{ahe1} and \eqref{ahe2} can be written respectively as
\beq\label{heq1}
b^{k}_{lj} c^{l}_{im}-b^{k}_{li} c^{l}_{jm}=a^k_{lj} c^{l}_{im}-a^k_{li} c^{l}_{jm},
\eeq
and
\beq\label{heq2}
b^{k}_{lj} c^{*l}_{im}-b^{k}_{li} c^{*l}_{jm}=a^k_{lj} c^{*l}_{im}-a^k_{li} c^{*l}_{jm}.
\eeq
This implies
$$b^{k}_{lj} c^{*l}_{im}E^i-b^{k}_{li} c^{*l}_{jm}E^i=a^k_{lj} c^{*l}_{im}E^i-a^k_{li} c^{*l}_{jm}E^i.$$
Taking into account that $E$ is the unity of the dual product and that $\nabla^* E=0$ one obtains
\beq\label{dualfromnatural}
b^{k}_{ij} =a^k_{ij}- c^{*l}_{ji}\nabla_l E^k \quad \text{   or   } \quad \nabla^*_XY=\nabla_X Y-\nabla_{(X*Y)}E,
\eeq
and similarly from \eqref{ahe1} one gets
\beq\label{naturalfromdual}
a^{k}_{ij} =b^{k}_{ij}- c^{l}_{ji}\nabla^*_l e^k, \quad \text{   or   } \quad \nabla_XY=\nabla^*_XY-\nabla_{(X\circ Y)}e.
\eeq
Using the definition of $*$
and applying to both sides the morphism $L:=E\circ: \T \rightarrow \T$ or replacing $Y$ with $Y'=E\circ Y$, one obtains the following in terms of the structure functions and the components of $L$:
\begin{equation}\label{aux10.eq}c^i_{lk}=L^i_jc^{*j}_{lk}\qquad\text{or}\qquad c^i_{lk}=L^j_lc^{*i}_{jk},\end{equation}
or equivalently 
\begin{equation}\label{aux10noindex.eq} 
X\circ Y=L(X *Y) \quad \text{ or } \quad X\circ Y=(LX*Y).
\end{equation}
\begin{remark}\label{remark1}
The first axiom in the definition of of a flat F-manifold is equivalent to the following two conditions: the flatness of $\nabla$ and the equation:
\begin{equation}\label{nablac.eq}
Y(\na_X\circ) Z=X(\na_Y\circ) Z\quad \text{ or } \quad \nabla_i c^l_{kj}=\nabla_kc^l_{ji}.
\end{equation}
Similarly, by definition of bi-flat F-manifold: 
\begin{equation}
Y(\nabla^*_X *)Z=X(\nabla^*_Y *)Z \quad \text{ or } \quad \nabla^*_i c^{*l}_{kj}=\nabla^*_kc^{l*}_{ji}.
\end{equation}
\end{remark}
We introduce formally Dubrovin-Frobenius manifolds without reference to flat F-manifolds and then we comment how they can be thought as special classes of flat (bi-flat really) F-manifolds. 
\begin{definition}\label{DFdefi}
A Dubrovin-Frobenius manifold $(M, \circ, \eta, e, E)$ is a manifold $M$ equipped with $\O$-bilinear sheaf morphism $\circ: \T \times \T \rightarrow \T$, two distinguished sections $e, E$ of $\T$ (called respectively unit and Euler vector field) and a flat non-degenerate pseudo-Riemannian metric $\eta$ (in the holomorphic case this is just a complex valued non-degenerate symmetric bilinear form) such that:
\begin{itemize}
\item the sheaf morphism $\circ$ induces a commutative associative product on $\T$;
\item the bilinear form $\eta$ is invariant with respect to $\circ$, i.e. $\eta(X\circ Y, Z)=\eta(X, Y\circ Z),$ for all local sections $X,Y,Z$ of $\T$;
\item indicating with $\nabla$ the Levi-Civita connection of $\eta$, $\nabla +\lambda \circ$ is a one-parameter family of connections which is flat and torsionless for any $\lambda$.;
\item for every local section $X\in \T$, $e\circ X=X$, and $\nabla e=0$;
\item $[E, e]=E$, $\mathcal{L}_E \circ =\circ,$ and $\mathcal{L}_E \eta=D\eta$, where $\mathcal{L}_E$ is the Lie derivative  and $D$ is a constant.
\end{itemize}
\end{definition}

At the points  where the operator $L=E\circ$ is invertible, any Dubrovin-Frobenius manifold can be endowed with a second flat contravariant metric, called \emph{the intersection form}, obtained applying $L$  to the inverse of the invariant metric $\eta$. This leads to the notion of almost duality for a Frobenius manifold (\cite{Dad}):
\begin{theorem}\label{DFAlmostDualdefi}
Let $(M, \circ, e, E, \eta)$ be a Dubrovin-Frobenius manifold. On the locus where $E\circ: \T\rightarrow \T$ is an isomorphism of sheaves, it is possible to induce a new product on local sections of $\T$ via the formula
\[X*Y=(E\circ)^{-1}X\circ Y.  \] 
Furthemore, it is possible to introduce a new metric (called intersection form) on local sections of $\T^*$ via:
\[ g(\alpha, \beta) =(\eta(\alpha) \circ E)(\beta), \quad \alpha, \beta \in \T^*,\]
where $\eta(\alpha)\in \T$ is the local vector field obtained by evaluating the (contravariant) metric $\eta$ on the local form $\alpha$. 
The product induced by $*$ is commutative and associative, with unit given by $E$. Furthermore, $g$ (covariant metric) is invariant with respect to the dual product, i.e. $g(X*Y, Z)=g(X, Y*Z)$ for all $X, Y, Z\in \T$ and the one-parameter family of connections $\nabla^*+\lambda*$ is torsionless and flat for any $\lambda$, where $\nas$ is the Levi-Civita connection of $g$. 
\end{theorem}
It turns out that the almost dual structure of a Dubrovin-Frobenius manifold, is a Dubrovin-Frobenius manifold without Euler vector field but the condition $\nabla^*E=0$ is not satisfied in general, see \cite{Dad}.

Given a Dubrovin-Frobenius manifold, on the locus where the almost dual structure exists, it is defined a so called {\em pencil of flat (contravariant) metrics} $g^{ij}_{\lambda}:=g^{ij}-\lambda \eta^{ij}$ (see \cite{du97}), i.e. this is a flat contravariant metric for any $\lambda$, and the contravariant Christoffel symbols $\Gamma_{\lambda,\, i}^{jk}:=-g_{\lambda}^{jl}\Gamma_{(\lambda)\, li}^k$ of $g_{\lambda}$ (where $\Gamma_{(\lambda)\, li}^k$ are the usual Christoffel symbols of $g_{\lambda,\, ij}$) are the linear combination of the contravariant Christoffel symbols of $g$ and $\eta$ respectively, i.e. $\Gamma_{\lambda i}^{jk}=\Gamma_{(g)\, i}^{jk}-\lambda \Gamma_{(\eta)\, i}^{jk}$. 

Given these structures, one can look for flat coordinates of $g_{\lambda}$ (these are the so called $\lambda$-periods). These are given by the extended Gauss-Manin system 
\[(g^{hj}-\lambda\eta^{hj})\d_j(d\theta)_k+(\Gamma^{hj}_{(g) \, k}-\lambda \Gamma^{hj}_{(\eta) \, k})(d\theta)_j=0.\]
It turns out that for Dubrovin-Frobenius manifold, $L^i_j=g^{il}\eta_{lj}$, where $L=E\circ$. Using this relation, the system for the $\lambda$-periods can be also written as
\begin{equation}\label{GMDF.eq}
(L-\lambda I)^j_i\d_j(d\theta)_k-(L^{l}_i\Gamma^{j}_{(g) \, lk}-\lambda \Gamma^{j}_{(\eta) \,ik})(d\theta)_j=0
\end{equation}
or (multiplying by $(L-\lambda I)^{-1}$) as
\[\d_h(d\theta)_k-\Gamma^j_{(g), hk}(d\theta)_j-\lambda((L-\lambda I)^{-1})^s_h(\Gamma^j_{(g)\, sk}- \Gamma^{j}_{\eta\, sk})(d\theta)_j=0.\]
The above system can be interpreted as the system for the flat coordinates for the connection defined by the Christoffel symbols:
\begin{equation}\label{GMDF2.eq}\Gamma^j_{(\lambda), \,hk}(\lambda)=\Gamma^j_{(g)\, hk}+\lambda(L-\lambda I)^{-1})^s_h(\Gamma^j_{(g)sk}- \Gamma^{j}_{(\eta)\, sk}).\end{equation}

 \section{A family of Gauss-Manin connections}\label{GM.section}
Motivated by the considerations above, we introduce the family of Gauss-Manin connections associated to any bi-flat F-manifold: 
\begin{definition}
Let $(M,\nabla,\circ,e,\nabla^*,*,E)$ be a bi-flat F-manifold.  We call \emph{Gauss-Manin connections}   the one-parameter family of connections defined by
\begin{equation}\label{GMnoindex}
\nabla^{GM}_XY=\nabla^*_XY+\lambda(\nabla^*_{L_{\lambda}^{-1}X}Y-\nabla_{L_{\lambda}^{-1}X}Y).
\end{equation}
where $L_{\lambda}=L-\lambda I, L=E\circ$. 
\end{definition}
Denoting by $a^i_{jk}$ and $b^i_{jk}$ the Christoffel symbols of the connections $\nabla$ and $\nabla^*$ respectively, the Christoffel symbols of \eqref{GMnoindex} can be written as
\begin{equation}\label{GM}
\Gamma^j_{hk}:=b^j_{hk}+\lambda((L_{\lambda})^{-1})^s_h(b^j_{sk}-a^{j}_{sk}),
\end{equation}
or, using \eqref{naturalfromdual},
\begin{equation}\label{GM2}
\Gamma^j_{hk}=b^j_{hk}+\lambda((L_{\lambda})^{-1})^s_hc^{l}_{sk}\nabla^*_le^j,
\end{equation}
or, using \eqref{dualfromnatural}
\begin{equation}\label{GM3}
\Gamma^j_{hk}=b^j_{hk}-\lambda((L_{\lambda})^{-1})^s_hc^{*l}_{sk}\nabla_l E^j,
\end{equation}

\begin{theorem}\label{GMth}
The  family of connections \eqref{GMtildenoindex} is flat and torsionless for any fixed $\lambda$ on the open set where $L_\lambda$ is invertible.  
\end{theorem}

\emph{Proof}. 
 It is enough to show that $T_{\nabla^{GM}} (\Ll X, \Ll Y)=0$ identically, where $T_{\nabla^{GM}}$ is the torsion of $\nabla^{GM}$. We have:
 \begin{eqnarray*}
 T_{\nabla^{GM}}(\Ll X, \Ll Y)&=&\nabla^{GM}_{\Ll X} \Ll Y-\na^{GM}_{\Ll Y}\Ll X -[\Ll X, \Ll Y]\\
 &\stackrel{\eqref{dualfromnatural}}{=}&T_{\nabla^*}(\Ll X, \Ll Y)-\lambda(\nabla_{X*\Ll Y} E)+\lambda (\nabla_{Y*\Ll X} E)\\
 &=& \lambda(\na_{Y* LX} E-\na_{X* LY} E)+\lambda^2(\na_{X* Y} E-\na_{Y* X} E),
 \end{eqnarray*}
 which vanishes identically because $T_{\nabla^*}=0$, $*$ is commutative, \eqref{aux10noindex.eq} and the fact that  $\circ$ is also commutative. 

To prove that the curvature $R^{GM}$ of $\nabla^{GM}$ is zero, we use formula \eqref{GMnoindex}. Since we work on the open dense subset where $L_{\lambda}$ is invertible, it is enough to compute
\[R^{GM}(\Ll X,\Ll Y)Z=\na^{GM}_{\Ll X}\na^{GM}_{\Ll Y}Z-\na^{GM}_{\Ll Y}\na^{GM}_{\Ll X}Z-\nabla^{GM}_{[\Ll X, \Ll Y]}Z,\]
We have 
\begin{eqnarray*}
\nabla^{GM}_{\Ll X}\na^{GM}_{\Ll Y}Z&=&\na^{GM}_{\Ll X}\left(\nas_{\Ll Y}Z+\lambda(\nas_Y Z-\na_Y Z) \right)\\
&=&\nas_{\Ll X}\nas_{\Ll Y}Z+\lambda \nas_{\Ll X}(\nas_Y Z-\na_Y Z)\\
&&+\lambda(\nas_X \nas_{\Ll Y} Z-\na_X \nas_{\Ll Y} Z)\\
&&+\lambda^2\left[ \nas_X(\nas_YZ -\na_Y Z)-\na_X(\nas_Y Z-\nabla_YZ) \right].
\end{eqnarray*}
On the other hand, the term $\nabla^{GM}_{[\Ll X, \Ll Y]}Z$ gives:
\begin{eqnarray*}
\nabla^{GM}_{[\Ll X, \Ll Y]}&=&\nas_{[\Ll X,\Ll Y]}Z+\lambda(\nas_{\Ll^{-1}[\Ll X, \Ll Y]} Z-\na_{\Ll^{-1}[\Ll X, \Ll Y]} Z)=\\
&&\lambda\left[ \nas_{\Ll^{-1} N_{L}(X, Y)} Z+\nas_{[LX, Y]+[X, LY]-L[X, Y]}Z-\lambda\nas_{[X,Y]} Z \right]\\
&&-\lambda\left[ \na_{\Ll^{-1} N_{L}(X, Y)}Z+\na_{[LX, Y]+[X, LY]-L[X, Y]}Z-\lambda\na_{[X,Y]} Z \right]\\
&&+\nas_{[\Ll X, \Ll Y]}Z,
\end{eqnarray*}
where $N_{L}$ is the Nijenhuis torsion of $L$, we have used the fact that the $N_{L}=N_{\Ll}$ and we expanded $\Ll=L-\lambda \mathbb{I}$ and the corresponding commutators. Using the Riemann tensor $R_{\nabla}$ of $\nabla$ and the Riemann tensor $R_{\nabla^*}$ of $\nabla^*$, we obtain 
\begin{eqnarray*}
&&R^{GM}(\Ll X, \Ll Y)Z=R_{\nabla^*}(\Ll X, \Ll Y)Z\\
&&+\lambda(R_{\nabla^*}(\Ll X,Y)Z+\nabla^*_{[\Ll X,Y]}Z)+\lambda(R_{\nabla^*}(X,\Ll Y)Z+\nabla^*_{[X,\Ll Y]}Z)\\
&&-\lambda\left[ \nas_{\Ll^{-1} N_{L}(X, Y)} Z+\nas_{[LX, Y]+[X, LY]-L[X, Y]}Z-\lambda\nas_{[X,Y]} Z \right]\\
&&+\lambda\left[ \na_{\Ll^{-1} N_{L}(X, Y)}Z+\na_{[LX, Y]+[X, LY]-L[X, Y]}Z-\lambda\na_{[X,Y]} Z \right]\\
&&+\lambda(-\nabla^*_{L_{\lambda}X}\nabla_YZ-\nabla_X\nabla^*_{L_{\lambda}Y}Z+\nabla^*_{L_{\lambda}Y}\nabla_XZ+\nabla_Y\nabla^*_{L_{\lambda}X}Z)\\
&&+\lambda^2(R_{\nabla^*}(X,Y)Z+\nabla^*_{[X,Y]}Z+R_{\nabla}(X,Y)Z+\nabla_{[X,Y]}Z)\\
&&+\lambda^2(-\nabla^*_X\nabla_YZ-\nabla_X\nabla^*_YZ+\nabla^*_Y\nabla_XZ+\nabla_Y\nabla^*_XZ).
\end{eqnarray*}
Expanding again $\Ll=L-\lambda\mathbb{I}$ and using the vanishing  of  the curvature
 of $\nabla$ and $\nabla^*$ we get first
\begin{eqnarray*}
&&R^{GM}(\Ll X, \Ll Y)Z=\\
&&+\lambda\left[\nas_{L[X, Y]}Z+\na_{[LX, Y]+[X, LY]-L[X, Y]}Z\right]\\
&&+\lambda(-\nabla^*_{L_{\lambda}X}\nabla_YZ-\nabla_X\nabla^*_{L_{\lambda}Y}Z+\nabla^*_{L_{\lambda}Y}\nabla_XZ+\nabla_Y\nabla^*_{L_{\lambda}X}Z)\\
&&+\lambda^2(-\nabla^*_X\nabla_YZ-\nabla_X\nabla^*_YZ+\nabla^*_Y\nabla_XZ+\nabla_Y\nabla^*_XZ),
\end{eqnarray*}
and then, after some additional computations,
\begin{eqnarray}
\label{GMR}
&&R^{GM}(\Ll X, \Ll Y)Z=\\
\nonumber&&\lambda (\nas_{L Y}\na_X Z-\na_X \nas_{L Y} Z+\na_Y \nas_{L X} Z-\nas_{L X}\na_Y Z)\\
\nonumber&&+\lambda( \nas_{L[X,Y]}Z+\na_{[LX,Y]+[X,LY]-L[X,Y]}Z).
\end{eqnarray}
Now using relation \eqref{dualfromnatural}, we eliminate $\nas$ in the first two lines of the previous expression and we get:
\begin{eqnarray*}
&&R^{GM}(\Ll X, \Ll Y)Z=\\
&&\lambda\left[ \bcancel{-\na_X\na_{LY}X}+\na_X\na_{LY*Z}E-\cancel{\na_{LX}\na_YZ}+\na_{LX*\nabla_YZ}E\right]\\
&&+\lambda\left[ \cancel{\na_Y\na_{LX}Z}-\na_Y\na_{LX*Z}E+\bcancel{\na_{LY}\na_XZ}-\na_{LY*\na_XZ}E\right]\\
&&+\lambda\left[ \nas_{L[X,Y]}Z+\cancel{\na_{[LX,Y]}Z}+\bcancel{\na_{[X,LY]}Z}-\na_{L[X,Y]}Z\right],
\end{eqnarray*}
where again we have canceled terms that sum to zero due to the vanishing of the curvature of $\na$. 
Now, by equation \eqref{dualfromnatural} 
\[\nas_{L[X,Y]}Z-\na_{L[X,Y]}Z=-\na_{(L[X,Y])*Z} E=-\na_{[X,Y]\circ Z}E\] using \eqref{aux10noindex.eq}. Likewise, we convert all the $*$ products in $\circ$ products eliminating all the occurrences of $L$ in the remaining terms using again  \eqref{aux10noindex.eq}.
We are left with 
\begin{eqnarray*}
R^{GM}(\Ll X, \Ll Y)Z&=&\lambda \na_{X\circ \na_Y Z}E-\lambda\na_Y \na_{X\circ Z}E-\lambda\na_{Y\circ \na_X Z}E\\
&&+\lambda\na_X \na_{Y\circ Z}E-\lambda\na_{[X,Y]\circ Z}E.
\end{eqnarray*} 
Using the second covariant derivative formula
\begin{equation}\label{scov.eq}
\nabla^2_{X,Y}Z=\nabla_X \nabla_Y -\nabla_{\nabla_X Y} Z
\end{equation}
and taking into account that $\nabla^2 E=0$, we have
\begin{equation}\label{aux12.eq}
\nabla_Y \nabla_{X\circ Z} E=\nabla_{\na_Y (X\circ Z)}E, \quad \nabla_X\nabla_{Y\circ Z} E=\nabla_{\na_X (Y\circ Z)}E.
\end{equation}
Substituting these in the previous expression for $R^{GM}$ we have
\begin{eqnarray*}
R^{GM}(\Ll X, \Ll Y)Z&=&\lambda \na_{X\circ \na_Y Z}E-\lambda\nabla_{\na_Y (X\circ Z)}E-\lambda\na_{Y\circ \na_X Z}E\\
&&\lambda\nabla_{\na_X (Y\circ Z)}E-\lambda\na_{[X,Y]\circ Z}E,
\end{eqnarray*} 
which after a straightforward computation vanishes identically using equation \eqref{nablac.eq} and the fact that the torsion of $\nabla$ is zero. 

\endproof

\begin{remark}\label{alternativeproof}
Alternatively, using the identity (that follows from Equation \eqref{heq1})
\[b^{k}_{sj} L^{s}_{i}-b^{k}_{si} L^{s}_{j}=a^k_{sj} L^{s}_{i}-a^k_{si} L^{s}_{j},\]
one can prove the identity
\begin{eqnarray}
\label{RiemId}
&&R^j_{lik}(L_{\lambda})^l_h(L_{\lambda})^i_t=\\
\nonumber&&(R_{\nabla^*})^j_{ilk}(L_{\lambda})^l_h(L_{\lambda})^i_t+\lambda^2((R_{\nabla^*})^j_{htk}-(R_{\nabla})^j_{htk})\\
\nonumber&&+\lambda\left(L^l_h((R_{\nabla^*})^j_{tlk}-(R_{\nabla})^j_{tlk})+L^l_t((R_{\nabla^*})^j_{hlk}-(R_{\nabla})^j_{hlk})\right)\\
\nonumber&&+\lambda(\nabla_t(L^l_h\Delta^j_{lk})-\nabla_h(L^l_t\Delta^j_{lk}))-\lambda (N_L)^m_{th}\left((L_{\lambda})^{-1}\right)^s_m\Delta^j_{sk}.
\end{eqnarray}
where
\[\Delta^j_{lk}=b^j_{lk}-a^j_{lk}\stackrel{\eqref{dualfromnatural}}{=}-c^{*s}_{lk}\nabla_sE^j\stackrel{\eqref{aux10.eq}}{=}-(L^{-1})^r_lc^{s}_{rk}\nabla_sE^j.\]
This implies
\begin{eqnarray*}
R^j_{lik}(L_{\lambda})^l_h(L_{\lambda})^i_t&=&(R_{\nabla^*})^j_{ilk}(L_{\lambda})^l_h(L_{\lambda})^i_t+\lambda^2((R_{\nabla^*})^j_{htk}-(R_{\nabla})^j_{htk})\\
&&+\lambda\left(L^l_h((R_{\nabla^*})^j_{tlk}-(R_{\nabla})^j_{tlk})+L^l_t((R_{\nabla^*})^j_{hlk}-(R_{\nabla})^j_{hlk})\right)\\
&&-\lambda(\nabla_tc^{s}_{hk}-\nabla_hc^{s}_{tk})\nabla_sE^j-\lambda(c^{s}_{hk}\nabla_t\nabla_sE^j-c^{s}_{tk}\nabla_h\nabla_sE^j)\\
&&+\lambda  (N_L)^m_{th}\left((L_{\lambda})^{-1}\right)^s_m(L^{-1})^r_sc^{t}_{rk}\nabla_tE^j.
\end{eqnarray*}
Taking into account that $R_{\nabla^*}=0,R_{\nabla}=0,N_L=0, \nabla\nabla E=0$ and that the  tensor $\nabla  c$ is symmetric  one obtains the result.
\end{remark}

The following result provides a converse to the previous statement, in the sense that starting from some data for which the family of Gauss-Manin connections is flat and torsionless, one can reconstruct a bi-flat F-manifold. 
\begin{theorem}\label{converseth}
Let $(M, \circ, e)$ be an F-manifold with an Euler vector field $E$ and a torsionless connection $\nabla$ on $\T$ (not assumed to be flat) such that $\nabla e=0$ and such that $\nabla E: \T\rightarrow \T$ is non-degenerate. Construct a new product $*$ such that $X*Y=(E\circ)^{-1}X\circ Y$  and a new connection $\nabla^*$ using \eqref{dualfromnatural}. 
Then the Gauss-Manin connection associated to these data is flat and torsionless iff the data $(M, \circ, e, \nabla, *, E,  \nabla^*)$ define a bi-flat F-structure. 
\end{theorem}
\proof
One direction is already provided by the Theorem \ref{GMth}. For the other direction we argue as follows. Since $E$ is an Euler vector field then $L:=E\circ$ has  zero Nijenhuis torsion. Construct the product $*$ and the connection $\nabla^*$ as in the statement of the theorem, and build the corresponding Gauss-Manin connection and assume it is torsionless and flat. From the vanishing of the torsion of $\nabla^{GM}$ it follows immediately that also $\nabla^*$ is torsionless. 

Now we look at the expression of the curvature $R^{GM}$ given in coordinates in the  Remark \ref{alternativeproof}. Since the Nijenhuis torsion of $L$ is zero,  $R^{GM}$ is a polynomial in $\lambda$, so in order for it to be zero all its coefficients have to be zero. In particular one deduces immediately that $R_{\nabla^*}=0$, $R_{\nabla}=0$. Furthermore, the flatness of $\nabla^*$ is equivalent to $\nabla^2 E=0$. Finally, the only remaining term for $R^{GM}$ is $\lambda(\nabla_tc^{s}_{hk}-\nabla_hc^{s}_{tk})\nabla_sE^j$ whose vanishing, due to the non-degeneracy of $\nabla_s E^j$, gives the condition $Y\nabla_X(\circ)Z =X\nabla_Y(\circ )Z$ which together with $R_{\nabla}=0$ implies that $(M, e, \circ, \nabla)$ is a flat F-manifold (see Remark \ref{remark1}). 

It remains to prove that $(M, *, E, \nabla^*)$ is a flat F-manifold. By construction $*$ induces a commutative associative product on $\T$ with unit $E$. Moreover, by definition $\nabla^*E=0$. So one needs only to check that $\nabla^*$ is compatible with $*$, i.e. $X\nabla^*_Y(*)Z=Y\nabla^*_X(*)Z$. This property is  inherited from the same property of the natural flat structure (see \cite{ALmulti}, unnumbered equation below equation 1.8 and the comment below it).
\endproof

\vspace{.5cm}
The results obtained so far in this section provide the two-sided arrow on the right in the triangular diagram in the Introduction. 

\subsection{The case  of Dubrovin-Frobenius manifolds}
The family  above can be further extended introducing an auxiliary parameter $\mu$ and replacing $\nabla^*$ with $\tilde\nabla^*=\nabla^*+\mu *$. It is immediate to check that $\tilde\nabla^*$ is torsionless and flat. The extended Gauss-Manin connection becomes:
\begin{equation}\label{GMtildenoindex}
\tilde\nabla^{GM}_XY=\tilde\nabla^*_XY+\lambda(\tilde\nabla^*_{L_{\lambda}^{-1}X}Y-\nabla_{L_{\lambda}^{-1}X}Y).
\end{equation}
with  Christoffel symbols 
\begin{eqnarray}\label{GM3tilde}
\tilde\Gamma^j_{hk}&=&b^j_{hk}+\mu c^{*j}_{hk}+\lambda((L_{\lambda})^{-1})^s_h(b^j_{sk}+\mu c^{*j}_{sk}-a^{j}_{sk})\\
\nonumber&=&b^j_{hk}+\mu c^{*j}_{hk} -\lambda((L_{\lambda})^{-1})^s_hc^{*l}_{sk}(\nabla E-\mu I)^j_l.
\end{eqnarray}
Using $\nas_XY=\na_X Y-\na_{X*Y}E$, we get the relation 
\begin{equation}\label{aux111.eq}
\tilde\nabla^*_XY=\nabla_X Y-\nabla_{(X*Y)}E+\mu(X*Y).
\end{equation}
The Theorem   \ref{GMth} is still valid and can be proved in a similar way. For the torsion we get, taking into account that the torsion of $\tilde\nabla^*$ is zero and equation \eqref{aux111.eq}:
\begin{eqnarray*}
 &&T_{\tilde\nabla^{GM}}(\Ll X, \Ll Y)=\\
 && \lambda(-\na_{X*\Ll Y} E+\mu(X*\Ll Y)+\na_{Y*\Ll X}E-\mu(Y*\Ll X))=\\
 && \lambda(-\na_{X*LY}E+\lambda \na_{X*Y}E+\na_{Y*LX E}-\lambda \na_{Y*X}E)=0,
 \end{eqnarray*}
 where we used, $\Ll=L-\lambda\mathbb{I}$, the commutativity of $*$, the fact that $X*LY=X\circ Y$ and the commutativity of $\circ$. 
Moreover, taking into account that replacing $b^j_{sk}$ with  $b^j_{sk}+\mu c^{*j}_{sk}$ the condition \eqref{heq2} is still valid as it is immediate to see, and therefore also of \eqref{heq1} which is equivalent to it.  Thus for the  tensor field  $\tilde{R}^j_{lik}(L_{\lambda})^l_h(L_{\lambda})^i_t$ one obtains the same expression  \eqref{RiemId} with 
\begin{itemize}
\item $R_{\nabla^*}$  replaced by the Riemann tensor of the deformed connection $\tilde\nabla^*=\nabla^*+\mu *$, which also vanishes, since for a bi-flat F-manifold, the structure $(M, \nabla^*, *, E)$ in particular is a flat F-manifold;
\item the  tensor field $\Delta^j_{lk}$ replaced by the tensor field
\[-c^{*s}_{lk}(\nabla E-\mu I)^j_s.\]
\end{itemize}
Both changes do not affect the result, as a straightforward computation using \eqref{aux10.eq} and \eqref{nablac.eq} shows. 
\newline
\newline
In the case of Dubrovin-Frobenius manifolds the functions $a^k_{ij}$ coincide with the Christoffel  symbols $\Gamma^k_{(\eta), ij}$ of the Levi-Civita connection of the invariant metric $\eta$ and,   for a suitable  choice of $\mu=\bar\mu$, $b^k_{ij}+\bar{\mu}c^{*j}_{ij}$ are the Christoffel  symbols $\Gamma^k_{(g), ij}$ of the Levi-Civita connection of the intersection form (see for instance \cite{ALmulti}). Using  this fact it is immediate to check that
 for this choice of $\mu$ the Gauss-Manin connection \eqref{GMtildenoindex} coincides with the connection \eqref{GMDF2.eq}.
 
\section{A differential bicomplex associated with bi-flat structures} 
Consider the sheaf of forms $\OO^{\bullet}$ on the manifold $M$ with the usual grading. 
There are sheaf maps $d$, $d_L$: $\OO^k \rightarrow \OO^{k+1}$ defined on local sections as:
$$(d \omega)(X_0, \dots, X_k)=\sum_{i=0}^k (-1)^i X_i(\omega(X_0, \dots, \hat{X}_i, \dots, X_k))+$$
$$+\sum_{0\leq i<j\leq k}(-1)^{i+j}\omega([X_i, X_j], X_0, \dots, \hat{X}_i, \dots, \hat{X}_j, \dots X_k),$$
where $X_i(\omega(X_0, \dots, \hat{X}_i, \dots, X_k))$ denotes the action of the vector field $X_i$ on the function $\omega(X_0, \dots, \hat{X}_i, \dots, X_k)$,
 $$(d_L \omega)(X_0, \dots, X_k)=\sum_{i=0}^k (-1)^i (LX_i)(\omega(X_0, \dots, \hat{X}_i, \dots, X_k))+$$
$$+\sum_{0\leq i<j\leq k}(-1)^{i+j}\omega([X_i, X_j]_L, X_0, \dots, \hat{X}_i, \dots, \hat{X}_j, \dots X_k),$$
where $(LX_i)(\omega(X_0, \dots, \hat{X}_i, \dots, X_k))$ indicates the action of the vector field $LX_i$ obtained applying the endomorphism $L:\T \rightarrow \T $ to $X_i$ and 
\[[X_i,X_j]_L=[LX_i, X_j]+[X_i, LX_j]-L[X_i,X_j].\]
According to the theory of Fr\"{o}licher-Nijenhuis 
 \cite{FN}, if $L$ has vanishing Nijenhuis torsion then
\begin{eqnarray*}
d\cdot d_L+d_L\cdot d=0\quad \text{ and } \quad d_L^2=0.
\end{eqnarray*} 
The maps $d$ and $d_L$  can be extended to forms with values in the tangent sheaf $\T$, once a connection on $\T$ has been chosen. The extension of $d$, $\dna$ is the exterior covariant derivative we met in the last section. 
 Its main properties are (see for instance \cite{Lee}):
\begin{enumerate}
\item The operator $\dna$ coincides with $d$ when restricted to scalar valued forms. 
\item If $\omega \in \T\otimes \OO^0$ then $(\dna \omega)(X)=\na_X \omega$, where $X$ is any local section of $\T$. 
\item If the connection $\nabla$ is flat, then $\dna\circ \dna=0$ identically. 
\end{enumerate}

To define the extension of $d_L$, $\dlna$ we follow \cite{ALimrn}: given an endomorphism of the tangent sheaf $L\in\T^*\otimes_{\O}\T$ on a manifold $M$ endowed with a connection $\nabla$, we define the {\em $L$-exterior covariant derivative} $\dlna$ acting on $\T\otimes_{\O}\OO^k$ as follows
$$(\dlna \omega)(X_0, \dots, X_k)=\sum_{i=0}^k (-1)^i \na_{LX_i}(\omega(X_0, \dots, \hat{X}_i, \dots, X_k))+$$
$$+\sum_{0\leq i<j\leq k}(-1)^{i+j}\omega([X_i, X_j]_L, X_0, \dots, \hat{X}_i, \dots, \hat{X}_j, \dots X_k).$$
Its main properties are (see \cite{ALimrn} for details): 
\begin{enumerate}
\item $\dlna$ coincides with $d_L$ when restricted to scalar valued differential forms. 
\item If $\omega\in \T\otimes_{\O}\OO^0$, then $(\dlna \omega)(X)=\na_{LX} \omega$, where $X$ is any local section of $\T$. 
\item If the connection $\nabla$ is flat and if $L$ has zero Nijenhuis torsion, then $\dlna\circ \dlna=0$ identically. 
\end{enumerate}

Furthermore, the following also holds:
\begin{lemma}\label{aux1.prop}
Given a manifold $M$ with a connection $\nabla$ on $\T$ and $L\in \T^*\otimes \T$, the morphisms of sheaves $d_{\nabla}: \T\otimes_{\O}\OO^k \rightarrow \T\otimes_{\O}\OO^{k+1}$ and $d_{L\nabla}:\T\otimes_{\O}\OO^k\rightarrow \T\otimes_{\O}\OO^{k+1}$ introduced above, satisfy:
\begin{itemize}
\item For $\alpha\in \OO^k$ and $\beta\in\T\otimes_{\O}\OO^l$ we have
\begin{equation}\label{aux4.eq}
\begin{array}{r@{}l}
 d_{\nabla}(\alpha\wedge\beta) &{}=d\alpha\wedge\beta+(-1)^k\alpha\wedge d_{\nabla}\beta,\\
 d_{\nabla}(\beta\wedge\alpha) &{}=d_{\nabla}\beta\wedge \alpha+(-1)^l\beta\wedge d\alpha,
\end{array}
\end{equation}
and the same relations hold for $d_{L\nabla}$ in place of $d_{\nabla}$.
\item If $\alpha\in \OO^k$ and $\sigma\in \T\otimes_{\O}\OO^0\cong \T$, then 
\begin{equation}\label{aux5.eq} d_{\nabla}(\sigma\otimes \alpha)=(\nabla \sigma)\wedge \alpha+\sigma\otimes d\alpha,\end{equation}
and the same relation holds for $d_{L\nabla}$ in place of $d_{\nabla}$.
\end{itemize}
\end{lemma}
\proof 
For $d_{\nabla}$ see \cite{Lee}, Theorem 12.57. A similar proof also holds for $d_{L\nabla}$.
\endproof

Given a bi-flat structure $(\nabla,\circ,e,\nabla^* *,E)$ we consider the differential $d_{\nabla}$ associated with $\nabla$ and the differential $d_{L\nabla^*}$ associated with
 $\nabla^*$ and $L=E\circ$. For the definition of bicomplex see for instance \cite{W}. In our case, we do not have a natural bi-gradation of the sheaf of $\O$-modules  $\T\otimes_{\O}\OO^{\bullet}$, but we still call it a bicomplex following for instance \cite{DMH}. 

\begin{theorem}\label{BidifferentialTh}
On any bi-flat $F$-manifold $M$, the sheaf morphisms $d_{\nabla}$ and $d_{L\nabla^*}$ determine a differential bicomplex structure on the sheaf of graded $\O$-modules $\T\otimes_{\O}\OO^{\bullet}$.
\end{theorem}
\proof
We have that $d_{\nabla}^2$ vanishes identically on $\T\otimes_{\O}\OO^{\bullet}$, due to the vanishing of the curvature of $\nabla$ and the same holds for $d_{L\nabla^*}^2$, due to the vanishing of the curvature of $\nabla^*$ and the vanishing of the torsion of $L$ (see \cite{ALimrn}, Proposition 5.3 point (4)).
In order to prove that $(d_{\nabla}, d_{L\nabla^*})$ give rise to a differential bicomplex, one needs to show that $\mathcal{D}:=d_\nabla\circ d_{L\nabla^*}+d_{L\nabla^*}\circ d_\nabla$ vanishes identically on $\T\otimes_{\O}\OO^{\bullet}$. First we prove that it is enough to show that $\mathcal{D}$ vanishes identically on $\T\otimes_{\O}\OO^0\cong \T$. 
Observe that every local section of $\T\otimes_{\O}\OO^k$ can be written as a finite sum $\sum_i \sigma_i\otimes \omega_i$ where $\omega_i$s are local $k$-forms and $\sigma_i$s are local sections of $\T$. Thus it is enough to show that $\mathcal{D}(\sigma\otimes \omega)=0$ where $\omega$ is a local $k$-form and $\sigma$ a local vector field.
We compute 
\[
d_{\nabla}\cdot d_{L\nabla^*}(\sigma\otimes \omega)\stackrel{\eqref{aux5.eq}}{=}d_\nabla(d_{L\nabla^*}\sigma \wedge \omega+\sigma\otimes d_L\omega)\stackrel{\eqref{aux4.eq}}{=}\]
\[=(d_{\nabla}d_{L\nabla^*}\sigma)\wedge \omega+(-1)d_{L\nabla^*}(\sigma)\wedge d\omega+d_{\nabla}(\sigma)\wedge d_L\omega+\sigma\otimes dd_L\omega,
\]
\[
d_{L\nabla^*}\cdot d_\nabla(\sigma\otimes \omega)\stackrel{\eqref{aux5.eq}}{=}d_{L\nabla^*}(d_\nabla(\sigma)\wedge \omega+\sigma\otimes d\omega)\stackrel{\eqref{aux4.eq}}{=}\]
\[=d_{L\nabla^*}d_{\nabla}(\sigma)\wedge\omega+(-1)d_\nabla(\sigma)\wedge d_L\omega+d_{L\nabla^*}(\sigma)\wedge d\omega+\sigma\otimes d_Ld\omega.\]

Summing the two previous expressions we get that \[\mathcal{D}(\sigma\otimes \omega)=(d_{\nabla}\cdot d_{L\nabla^*}+d_{L\nabla^*}\cdot d_{\nabla})(\sigma)\wedge \omega+\sigma\otimes (d\cdot d_L+d_L\cdot d)(\omega)=\] \[=(d_{\nabla}\cdot d_{L\nabla^*}+d_{L\nabla^*}\cdot d_{\nabla})(\sigma)\wedge \omega,\] due to the vanishing of the torsion of $L$. 
Thus it is enough to show that $((d_\nabla\cdot d_{L\nabla^*}+d_{L\nabla^*}\cdot d_\nabla)(Z))^j_{ri}=0$ for each local section $Z$ of the tangent sheaf. 
We have
\begin{eqnarray*}
(d_\nabla \cdot d_{L\nabla^*} Z)^j_{ri}&=&\nabla_r(L^m_i\nabla^*_mZ^j)-\nabla_i(L^m_r\nabla^*_mZ^j)\\
&=&\partial_r(L^m_i\nabla^*_mZ^j)+a^j_{rl}L^m_i\nabla^*_mZ^l-\partial_i(L^m_r\nabla^*_mZ^j)-a^j_{il}L^m_r\nabla^*_mZ^l,\\
(d_{L\nabla^*}\cdot d_{\nabla} Z)^j_{ri}&=&L^m_r(\nabla^*_m\nabla_{(i)}Z^j)-L^m_i(\nabla^*_m\nabla_{(r)}Z^j)-\nabla_mZ^j(\partial_r L^m_i-\partial_iL^m_r)\\
&=&L^m_r(\nabla^*_m\nabla_{i}Z^j)-L^m_i(\nabla^*_m\nabla_{r}Z^j)-\nabla_mZ^j(\nabla^*_r L^m_i-\nabla^*_iL^m_r).
\end{eqnarray*}
The formulas above can be written as
\begin{eqnarray*}
(d_\nabla \cdot d_{L\nabla^*} Z)(X,Y)&=&\nabla_X\nabla^*_{LY}Z-\nabla_Y\nabla^*_{LX}Z-\nabla^*_{L[X,Y]},\\
(d_{L\nabla^*}\cdot d_{\nabla} Z)(X,Y)&=&\nabla^*_{LX}\nabla_YZ-\nabla^*_{LY}\nabla_XZ-\nabla_{\nabla^*_{LX}Y}Z\\
&&+\nabla_{\nabla^*_{LY}X}Z+\nabla_{\nabla^*_{Y}LX}Z-\nabla_{\nabla^*_{X}LY}Z+\nabla_{L[X,Y]}Z\\
&=&\nabla^*_{LX}\nabla_YZ-\nabla^*_{LY}\nabla_XZ+\nabla_{[Y,LX]+[LY,X]+L[X,Y]}Z.
\end{eqnarray*} 
Therefore we have
\begin{eqnarray*}
(d_\nabla \cdot d_{L\nabla^*}Z+d_{L\nabla^*}\cdot d_{\nabla} Z)(X,Y)&=&\nabla_X\nabla^*_{LY}Z-\nabla^*_{LY}\nabla_XZ\\
&&-\nabla_Y\nabla^*_{LX}Z+\nabla^*_{LX}\nabla_YZ\\
&&-\nabla^*_{L[X,Y]}+\nabla_{[Y,LX]+[LY,X]+L[X,Y]}Z
\end{eqnarray*}
and  this vanishes due to the vanishing of the curvature of the Gauss-Manin connections
 (compare with \eqref{GMR}).
\endproof

The following result provides, under suitable assumptions, a converse to the previous statement.
\begin{theorem}
Let $(M, \circ, e)$ be an F-manifold with an Euler vector field $E$ and a torsionless connection $\nabla$ on $\T$ (not assumed to be flat) such that $\nabla e=0$ and such that $\nabla E: \T\rightarrow \T$ is non-degenerate. Construct a new product $*$ such that $X*Y=(E\circ)^{-1}X\circ Y$  and a new connection $\nabla^*$ using \eqref{dualfromnatural}. Then the data $(d_{\nabla}, d_{L\nabla^*})$ give rise to a differential bicomplex on the sheaf of graded $\O$-modules $\T\otimes_{\O}\OO^{\bullet}$  iff the data $(\nabla, \circ, e, \nabla^*, *, E)$ define a bi-flat F-structure. 
\end{theorem}
\proof
One direction is already provided by the Theorem \ref{BidifferentialTh}. For the other direction, it is enough to follow the steps in the proof of Theorem \ref{converseth} and use the expression for $(d_\nabla \cdot d_{L\nabla^*}Z+d_{L\nabla^*}\cdot d_{\nabla} Z)(X,Y)$ at the end of Theorem \ref{BidifferentialTh} and compare it with the same expression that appear in the Remark  \ref{alternativeproof}.  We leave the details to the reader. 
\endproof

The results obtained so far in this section provide the two-sided arrow on the left in the triangular diagram in the Introduction. 

The differential bicomplex associated with a bi-flat F-manifold allows to define a sequence of vector fields by means of recursion relations that are similar to the classical  Lenard-Magri chain. Let $X_{(0)}$ be a local section of $\T$  for $\nabla$ satisfying the equation
\begin{equation}\label{cohoEq}
d_{\nabla}d_{L\nabla^*}X_{(0)}=0.
\end{equation}
Then we can define ``higher'' local sections of $\T$ by using the following recurrence relations
\beq\label{LMch}
d_{\nabla}X_{(\alpha+1)}=d_{L\nabla^*}X_{(\alpha)}, \quad \alpha\in \mathbb{N}
\eeq
or, equivalently
\begin{equation}\label{rec1.eq}  \na_Y X_{(\alpha+1)}= \nas_{LY} X_{(\alpha)}=\nas_{E\circ Y}X_{(\alpha)}.
\end{equation}

\begin{theorem}
The sequence $\{X_{(\alpha)}\}_{\alpha \in \mathbb{N}}$ of local sections of $\T$ is well defined.
\end{theorem}

The cohomology of $d_{\nabla}$ is trivial on a subset $U\subset M$ which is diffemorphic to a ball. Indeed, introducing flat coordinates for $\nabla$ in $U$ one has that $d_{\nabla}$ is simply the de Rham differential applied componentwise to a vector of forms. Thus the triviality of the cohomology of $d_{\nabla}$ follows from the triviality of the de Rham cohomology on a ball. Equation \eqref{cohoEq} tell us that $d_{L\nabla^*}X_{(0)}$
 is a $d_{\nabla}$ cocycle. Therefore there exists a vector field $X_{(1)}$ such that
\begin{equation*}\label{LMch0}
d_{\nabla}X_{(1)}=d_{L\nabla^*}X_{(0)}.
\end{equation*}
Since $d_{\nabla}$ and $d_{L\nabla^*}$ anticommute we have
\[d_{L\nabla^*}d_{\nabla}X_{(1)}=-d_{\nabla}d_{L\nabla^*}X_{(0)}=0.\]
Thus we can repeat the previous argument obtaining
\begin{equation*}\label{LMch1}
d_{\nabla}X_{(2)}=d_{L\nabla^*}X_{(1)},
\end{equation*}
and, in general, \eqref{LMch}.

\endproof

The sequence of vector field  defined recursively by \eqref{LMch}, starting from a frame $(X_{(1,0)},...,X_{(n,0)})$ of flat local sections of $\T$  for $\nabla$  (where $n=\text{dim}(M)$),  
 can be identified with the coefficients of 
 local flat section of the Gauss-Manin connections. 
Indeed, the following remarkable property for the local sections of $\T$ that are flat with respect to  $\nabla_{GM}$ holds (for simplicity  we assume $\mu=0$). It corresponds to the base arrow in the conceptual triangle places in the Introduction. 
\begin{theorem}\label{flatsecGM}
Let $X$ be a local section $X\in \T$ that is flat with respect to $\nabla_{GM}$. Then, viewing $X$ is as an element of $\T[[\lambda^{-1}]]$, then the formal power series coefficients satisfy the Lenard-Magri recurrence relations \eqref{LMch}. 
\end{theorem}
\proof

From $\nabla^{GM}X=0$, using arbitrary sections $Y\in \T$, we have
\begin{eqnarray*}
0=\nabla^{GM}_{\Ll Y}X&\stackrel{\eqref{GMnoindex}}{=}&\nabla^*_{\Ll Y}X+\lambda(\nabla^*_{Y}X-\nabla_Y X)\\
&=&\nabla^*_{L Y} X-\lambda\nabla^*_Y X+\lambda(\nabla^*_{Y}X-\nabla_Y X),\\
\end{eqnarray*}
which yields 
\begin{equation}\label{final.eq}
0=\nabla^{GM}_{\Ll Y}X=\nabla^*_{L Y} X-\lambda \nabla_Y X
\end{equation}
Now notice that $\T$ is made into a commutative associative ring thanks to $\circ$. Therefore it makes sense to consider $X\in\T$ such that $\nabla^{GM} X=0$ as an element of the ring of formal power series in $\lambda^{-1}$ with coefficients in the ring $\T$: $X\in \T[[\lambda^{-1}]]$. 
Representing $X_{(p)}$ such that $\nabla^{GM}X_{(p)}=0$ as an element of $\T[[\lambda^{-1}]]$, we get:  
$X_{(p)}=X_{(p,0)}+\frac{X_{(p,1)}}{\lambda}+\frac{X_{(p,2)}}{\lambda^2}+\cdots$
and substituting it in \eqref{final.eq}, we get recurrence relations on the coefficients of the expansion (due to the vanishing of  the curvature they are well defined).
 It is immediate to check that they coincide with the Lenard-Magri recurrence relations \eqref{LMch}.
\endproof

\begin{remark}
Theorem \ref{flatsecGM} begs the question of what is the relationship, in the semisimple case, between the local sections of $\T$, $\{X_{(p, \alpha)}\}_{\alpha \in \mathbb{N}}$ that satisfy the recurrence \eqref{LMch}, and the sequence of of vector fields $\{ Y_{(p, \alpha)}\}_{\alpha \in \mathbb{N}}$ that satisfy the twisted Lenard-Magri chains in \cite{ALimrn}. They are related in the following way: $Y_{(p,0)}=0, X_{(p,0)}=-Y_{(p,1)}$ and in general $X_{(p,\alpha)}=LY_{(p, \alpha)} -Y_{(p,\alpha+1)}$. It is not difficult to prove using induction that if the $\{X_{(p, \alpha)}\}_{\alpha \in \mathbb{N}}$ satisfy the recurrence \eqref{LMch}, then the $Y_{(p,\alpha)}$s satisfy the twisted Lenard-Magri chain and conversely. 

Using the results of \cite{ALimrn},  it is also easy to check that, in the semisimple case, the systems of PDEs 
\[u^i_{t_{(p, \alpha)}}=(d_{\nabla} Z_{(p, \alpha)})^i_j u^j_x\]  defined by the vector fields 
\[Z_{(p,\alpha)}=-\sum_{k=1}^{\alpha}L^{k-1}X_{(p,\alpha-k)}\] 
are compatilble (i.e. the flows commute).

Furthermore, it is easy to check that,  in the Dubrovin-Frobenius case, the coefficient of of a $\nabla^{GM}$-flat 1-form $\alpha$ expressed in terms of power series in $\lambda^{-1}$ fulfill a recurrence which coincides with the bi-Hamiltonian recurrence (see \cite{ALimrn}, Proposition 10.3).
\end{remark}

\end{document}